\documentclass[12pt,reqno,twoside]{amsart}
\usepackage{amsmath,amssymb,amsfonts,amsthm,indentfirst,setspace}
\usepackage{xcolor}\usepackage{ulem}
\usepackage{graphicx}

\newtheorem{theorem}{Theorem}[section]

\newtheorem{corollary}{Corollary}[section]

\newcommand{\be}{\begin{equation}}
\newcommand{\ee}{\end{equation}}
\newcommand{\bea}{\begin{eqnarray}}
\newcommand{\eea}{\end{eqnarray}}
\newcommand{\eeas}{\end{eqnarray*}}
\newcommand{\beas}{\begin{eqnarray*}}

\begin{document}

\title
[conformally flat generalized Ricci recurrent manifolds]
{A classification of conformally flat generalized Ricci recurrent pseudo-Riemannian manifolds}
\author{Tee-How Loo \and Avik De}
\address{T. H. Loo\\
Institute of Mathematical Sciences\\
Universiti Malaya\\
50603 Kuala Lumpur\\
Malaysia}
\email{looth@um.edu.my}
\address{A. De\\
Department of Mathematical and Actuarial Sciences\\
Universiti Tunku Abdul Rahman\\
Jalan Sungai Long\\
43000 Cheras\\
Malaysia}
\email{de.math@gmail.com}

\footnotetext{A.D. and L.T.H. are supported by the grant FRGS/1/2019/STG06/UM/02/6.}

\begin{abstract}
Conformally flat pseudo-Riemannian manifolds with generalized Ricci recurrent, $(GR)_n$ structure are completely classified in this short report. A conformally flat generalized Ricci recurrent pseudo-Riemannian manifold is shown to be either a de Sitter space or an anti-de Sitter space. 
In particular, a conformally flat generalized Ricci recurrent spacetime must be either a de Sitter spacetime or an anti-de Sitter spacetime. 
\end{abstract}
\date{}
\maketitle

\section{\textbf{Introduction}}
The importance of the Generalized Ricci recurrent structure and its interaction with several gravity theories, like the standard theory of gravity \cite{avik}, the modified $f(R)$-theory \cite{avik-2}, the modified Gauss-Bonnet $f(R,G)$-theory \cite{frg} are well established. And of course, mathematically the structure raised curiosity among researchers in several occasions, details can be seen in \cite{avik} and the references therein. 

De et al. \cite{de} introduced the notion of $(GR)_n$ as an $n$-dimensional non-flat pseudo-Riemannian manifold whose Ricci tensor $R_{ij}$ satisfies the following:
\be \nabla _{i}R_{jl}= A_iR_{jl}+B_ig_{jl},\label{eqn:1}\ee
where $A_i$ and $B_i$ are two non-zero 1-forms. The structure is considered to be a generalization of Patterson introduced Ricci recurrent manifolds $R_n$ \cite{patterson} in which the Ricci tensor satisfies $\nabla_{i}R_{jl}= A_iR_{jl}$, with a non-zero 1-form $A_i$. Obviously, if the one-form $B_i$ vanishes, it reduces to a $R_n$. 

A $(GR)_4$ spacetime is a generalized Robertson Walker spacetime with Einstein fibre for a Codazzi type $R_{ij}$ \cite{pinaki}. In \cite{avik}, a conformally flat $(GR)_4$ was shown to be a perfect fluid. Its interaction with general relativistic cosmology was discussed thoroughly in the same paper. Continuing to this study, in \cite{avik-2}, the authors studied a conformally flat $(GR)_4$ with constant $R$ as a solution of $f(R)$-gravity theory. The presently accepted homogeneous and isotropic model of our universe, the Robertson-Walker spacetime is $(GR)_4$ if and only if it is Ricci symmetric. The equation of state (EoS) was shown to have $\omega = -1$. Several energy conditions were also analyzed and validated by current observational dataset. Very recently, the impact of $(GR)_4$ structure is investigated in modified Gauss-Bonnet, $f(R,G)$ theory of gravity \cite{frg}. The obtained results were examined for two particular $f(R,G)$-models and for them both, the weak, null and dominant energy conditions were validated while the strong energy condition was violated, which is a good agreement with the recent observational studies revealing that the current universe is in accelerating phase. Despite having such recognised physical significance, a complete classification of conformally flat Lorentzian manifolds with generalized Ricci recurrent structure is still due. 
This motivated us to do a careful scrutiny through all the possible cases.
We shall give a complete classification of this structure in a general setting as below:
\begin{theorem}\label{res1}
A conformally flat generalized Ricci recurrent pseudo-Riemannian manifold is an Einstein manifold, in particular, 
it is either a de Sitter space or an anti-de Sitter space.
\end{theorem} 
We can then deduce a classification of conformally flat generalized Ricci recurrent spacetime.
\begin{corollary}\label{res2}
A conformally flat generalized Ricci recurrent spacetime is either a de Sitter spacetime or an anti-de Sitter spacetime.
\end{corollary} 


\section{The Proof}
Let $M$ be an $n$-dimensional pseudo-Riemannian manifold, $n\geq 3$, with index $q$ ($0\leq q\leq n$).
Suppose $M$ is a conformally flat $(GR)_n$, that is, the Ricci tensor satisfies (\ref{eqn:1})
and the Riemannian curvature tensor satisfies \cite[pp.302]{neil}
\begin{align}\label{eqn:2}
R^l{}_{ijk}=\frac{R_{ki}\delta^l_j-R_{ji}\delta^l_k+g_{ki}R^l_j-g_{ji}R^l_k}{n-2}
-\frac{g_{ki}\delta^l_j-g_{ji}\delta^l_k}{n-1}\frac{R}{n-2}.
\end{align}
The divergence of the above equation gives
\begin{align}\label{eqn:3}
\nabla_jR_{ik}-\nabla_kR_{ij}=\frac12\frac{\nabla_jR g_{ki}-\nabla_kR g_{ji}}{n-1}.
\end{align}
Contracting $j$ and $l$ in (\ref{eqn:1}) we obtain
\be \nabla_iR=RA_i+nB_i.\label{eqn:6}\ee
Contracting  $i$ and $l$ in (\ref{eqn:1}) we obtain
\be \frac{1}{2}\nabla_jR=A^lR_{lj}+B_j.\label{eqn:7}\ee
It follows from (\ref{eqn:6})--(\ref{eqn:7}) that 
\be\label{eqn:8} 2A^lR_{lj}=RA_j+(n-2)B_j.\ee
By using (\ref{eqn:1}), (\ref{eqn:3})--(\ref{eqn:6}) we calculate
\be \label{eqn:9}
2(n-1)\{A_iR_{kl}-A_kR_{il}\}=RA_ig_{kl}-RA_kg_{il}-(n-2)B_ig_{kl}+(n-2)B_kg_{il}.
\ee
Next we split into two cases: $A^i$ is not a lightlike vector and $A^i$ is a lightlike vector.

\medskip\textbf{(I) $A^i$ is either a timelike or spacelike vector.}
Write 
\[
A_i=\varepsilon\alpha U_i; \quad U_lU^l=\varepsilon=\pm 1; \quad \alpha\neq0. 
\]
Then (\ref{eqn:8})--(\ref{eqn:9}) become
\begin{align}
 2\varepsilon\alpha U^lR_{ij}=&\varepsilon\alpha RU_j+(n-2)B_j,   \label{eqn:8b}\\
2(n-1)\varepsilon\alpha\{ U_iR_{kl}-U_kR_{il}\}
=&\varepsilon\alpha\{ RU_ig_{kl}- RU_kg_{il}\}-(n-2)B_ig_{kl}+(n-2)B_kg_{il}.  \label{eqn:9b}
\end{align}
Transvecting (\ref{eqn:9b}) with $U^i$ we obtain
\bea
2(n-1)\alpha R_{kl}
=\{\alpha R-(n-2)\beta\} g_{kl}+(n-2)\{\varepsilon\alpha RU_kU_l+U_lB_k+(n-1)U_kB_l\},
\label{eqn:10}
\eea
where $\beta=B_lU^l$.
The skew-symmetry part gives 
\[ 
U_lB_k+(n-1)U_kB_l   = U_kB_l+(n-1)U_lB_k,
\]
which implies that 
\be B_i={\varepsilon}\beta U_i, \quad (\beta=  U^lB_l).\label{eqn:11}\ee
Armed with this significant result, using (\ref{eqn:6}), (\ref{eqn:8b}) and (\ref{eqn:10}) we can conclude that the Ricci curvature tensor and the Ricci scalar satisfy
\begin{align} 
2(n-1)R_{kl}=&\{R-(n-2)\lambda\}g_{kl} +(n-2)\{R+n\lambda\}\varepsilon U_kU_l, \label{eqn:12}	\\
\nabla_iR=&\{R+n\lambda\}\alpha \varepsilon U_i, \label{eqn:13} 	\\
2U^lR_{lj}=&\{R+(n-2)\lambda\}U_j, \label{eqn:14}  
\end{align}
where 
\[
\lambda=\frac\beta\alpha.
\]
Covariantly differentiating (\ref{eqn:12}) with respect to $i$ we obtain
\begin{align*}
2(n-1)\nabla_iR_{kl}
=& \{\nabla_iR-(n-2)\nabla_i\lambda\}g_{kl}+(n-2)(\nabla_iR+n\nabla_i\lambda)\varepsilon U_kU_l \notag\\
			&	 +(n-2)(R+n\lambda){\varepsilon}\nabla_i(U_kU_l) \notag\\
=&(n-2)\{-\nabla_i\lambda g_{kl}+n\nabla_i\lambda\varepsilon U_kU_l+(R+n\lambda){\varepsilon}\nabla_i(U_kU_l)\} \notag\\
 &+\nabla_iR \{g_{kl}+(n-2)
					\varepsilon U_kU_l\}.
\end{align*} 
On the other hand,  (\ref{eqn:1}), (\ref{eqn:11})--(\ref{eqn:12}) give
\begin{align*}
2(n-1)\nabla_iR_{kl}
=&2(n-1)\alpha\varepsilon U_i R_{kl}+2(n-1)\lambda\alpha\varepsilon U_i g_{kl} \notag\\
=&\{R-(n-2)\lambda\}\alpha\varepsilon U_ig_{kl} +(n-2)\{R+n\lambda\}\alpha\varepsilon U_i \varepsilon U_kU_l \notag \\
	&+2(n-1)\lambda\alpha\varepsilon U_i g_{kl} \notag \\
=& (R+n\lambda)\alpha\varepsilon U_i	\{g_{kl}+(n-2)\varepsilon U_kU_l\} \notag\\
=&\nabla_iR \{g_{kl}+(n-2)\varepsilon U_kU_l\}.	
\end{align*} 
These two equations give
\begin{align*}
 -\nabla_i\lambda g_{kl}+n\nabla_i\lambda\varepsilon U_kU_l+(R+n\lambda){\varepsilon}\nabla_i(U_kU_l)=0.
\end{align*}
Transvecting with $U^k$ and $U^l$ we obtain
\begin{align}\label{eqn:22}
\nabla_i\lambda=0.
\end{align}
It follows from these two equations that 
\begin{align*}
(R+n \lambda)\nabla_i(U_kU_l)=0.
\end{align*}

We consider two subcases: $R+n\lambda\neq0$ and $R+n\lambda=0$.
\begin{description}
\item[a] If $R+n\lambda\neq0$, we have $\nabla_i U_l=0$ and so $U_lR_i^l=0$.
It follows from (\ref{eqn:14}) that 
\[
 R+(n-2)\lambda=0.
\]
It follows from this equation and (\ref{eqn:22}) that 
\[
\nabla_iR=-(n-2)\nabla_i\lambda=0.
\]
After applying this to (\ref{eqn:13}) gives
$R+n\lambda=0$. Hence we have $R=\lambda=0$ and so $R_{lk}=0$; violating the hypothesis of being $(GR)_n$. 
Hence this case is impossible.
\item[b] Suppose that  
\[ R+n\lambda=0.
\]
Then, using (\ref{eqn:12}) we obtain
\[ 
R_{kl}=-\lambda g_{kl},
\]
meaning that it is an Einstein space.
\end{description}

\medskip
\textbf{(II) $A^i$ is a lightlike vector.}
Take another lightlike vector $K^i$ such that  
\[
A_lA^l=K_lK^l=0; \quad A_lK^l=-1. 
\]
Transvecting (\ref{eqn:9}) with $A^i$, with the help of (\ref{eqn:8}), we obtain
\be
-A^jB_j g_{kl}+ RA_kA_l+A_lB_k+(n-1)A_kB_l=0.
\label{eqn:30}
\ee
The skew-symmetry part gives 
\[ 
A_lB_k= A_kB_l,
\]
which implies that $A^kB_k=0$ and 
\be B_i=-\lambda A_i, \quad (\lambda=K^lB_l).\label{eqn:31}\ee
Substituting these into (\ref{eqn:6}) and  (\ref{eqn:30}) respectively give
\begin{align} 
\nabla_iR=R-n\lambda=0,
\label{eqn:32}\end{align}
and so
\be\nabla_i\lambda=0.
\label{eqn:32b}\ee
By using (\ref{eqn:31})--(\ref{eqn:32}), we can simply (\ref{eqn:9}) as 
\be\label{eqn:33}
A_iR_{kl}-A_kR_{il}=\lambda A_ig_{kl}-\lambda A_kg_{il}.
\ee
Transvecting with $K^i$ and $K^l$ gives
\[
K^lR_{lk}=-\tau A_k+\lambda K_k; \quad (\tau=K^iK^jR_{ij}). 
\]
Transvecting (\ref{eqn:33}) with $K^i$  and applying the above equation gives
\be\label{eqn:34}
R_{lk}=\lambda g_{lk}+\tau A_lA_k.
\ee
Suppose 
\be\label{eqn:tau} \tau\neq0. \ee
Applying (\ref{eqn:32}) and (\ref{eqn:34}), we can simplfy (\ref{eqn:2}) as 
\begin{align}\label{eqn:35}
R^l{}_{ijk}=\frac\lambda{n-1}\{g_{ki}\delta^l_j-g_{ji}\delta^l_k\}
			+\frac\tau{n-2}\{A_kA_i\delta^l_j-A_jA_i\delta^l_k+g_{ki}A^lA_j-g_{ji}A^lA_k\}.
\end{align}
On the other hand, by substituting (\ref{eqn:31}) and (\ref{eqn:34}) into  (\ref{eqn:1}) give
\be\label{eqn:36}
\nabla_iR_{jk}=\tau A_iA_jA_k.
\ee
It follows from (\ref{eqn:2}), (\ref{eqn:32})--(\ref{eqn:32b}) and (\ref{eqn:34}) that
\[
\nabla_iR_{hj}=\nabla_hR_{ij}=\nabla_h(\tau A_iA_j).
\]
Covariantly differentiating (\ref{eqn:36}) with respect to $h$ gives
\[
\nabla_h\nabla_iR_{jk}=\nabla_h(\tau A_iA_j)A_k+\tau A_iA_j\nabla_hA_k
=\nabla_hR_{ij}A_k+\tau A_iA_j\nabla_hA_k,
\]
and so
\begin{align}
\tau A_jA_i\nabla_hA_k-\tau A_jA_h\nabla_iA_k
=&
(\nabla_h\nabla_i-\nabla_i\nabla_h)R_{jk}
=-R^l{}_{jhi}R_{lk}-R^l{}_{khi}R_{lj} \notag	
\end{align}
Next, by applying (\ref{eqn:34}), the properties $R_{kjhi}+R_{jkhi}=0$ and (\ref{eqn:tau}) we have  
\begin{align*}
A_jA_i\nabla_hA_k-A_jA_h\nabla_iA_k=&-R^l{}_{jhi}A_lA_k  + R^l{}_{khi}A_lA_j.
\end{align*}
Furthermore, by using (\ref{eqn:35}) we obtain
\begin{align}\label{eqn:37}
A_jA_i\Big\{\nabla_hA_k-\frac\lambda{n-1}g_{hk}\Big\}-A_jA_h\Big\{\nabla_iA_k-&\frac\lambda{n-1}g_{ik}\Big\} \notag\\
=&\frac\lambda{n-1}\{-g_{ij}A_h+g_{hj}A_i\}A_k.
\end{align}
Transvecting with $K^i$ and $K^j$ we obtain
\begin{align*}
\nabla_hA_k-\frac\lambda{n-1}g_{hk}=-A_h\Omega_k-\frac\lambda{n-1}K_hA_k; \quad 
\left(\Omega_k=K^i\nabla_iA_k-\frac\lambda{n-1}K_k\right).
\end{align*}
Substituting this into (\ref{eqn:37}) gives
\be\notag
\frac\lambda{n-1}\{-A_jA_iK_h+A_jA_hK_i\}A_k=\frac\lambda{n-1}\{-g_{ij}A_h+g_{hj}A_i\}A_k.
\ee
Since $\lambda\neq0$, it is descended to 
\[
-A_jA_iK_h+A_jA_hK_i=-g_{ij}A_h+g_{hj}A_i.
\] 
Transvecting with $K^h$ and $g^{ij}$ we obtain $1=n-1$. This is impossible. Hence 
(\ref{eqn:tau}) is false or we must have $\tau=0$. It follows from 
(\ref{eqn:34}) that it is also an Einstein space in this case.

Hence we conclude that $M$ is an Einstein space in both cases.
 As a conformaly flat Einstein space is of constant sectional curvature and the manifold is non-flat, 
We conclude that a conformally flat generalized Ricci recurrent spacetime is either a de Sitter space or an anti-de Sitter space.



\end{document}